# SPECTRAL CONDITIONS OF PANCYCLICITY FOR $t$-TOUGH GRAPHS

**Vladimir I. Benediktovich**
Institute of Mathematics of NASB
Surganov str. 11, 220072 Minsk, Belarus
vbened@im.bas-net.by

## ABSTRACT

More than 40 years ago Chvátal introduced a new graph invariant, which he called graph toughness. From then on a lot of research has been conducted, mainly related to the relationship between toughness conditions and the existence of cyclic structures, in particular, determining whether the graph is Hamiltonian and pancyclic. A pancyclic graph is certainly Hamiltonian, but not conversely. Bondy in 1976, however, suggested the "metaconjecture" that almost any nontrivial condition on a graph which implies that the graph is Hamiltonian also implies that the graph is pancyclic. We confirm the Bondy conjecture for $t$-tough graphs in the case when $t \in \{1; 2; 3\}$ in terms of the edge number, the spectral radius and the signless Laplacian spectral radius of the graph.

## 1 Introduction

Let $G = (V(G), E(G))$ be a finite simple (without loops or multiple edges) undirected $(n, m)$-graph with vertex set $V(G) = \{v_1, \ldots, v_n\}$ and edge set $E(G)$. The neighborhood of a vertex $v_i \in V(G)$ is denoted by $N_G(v_i)$, and its degree by $d_i = d_G(v_i) = |N_G(v_i)|$. Denote by $\delta(G)$ ($\Delta(G)$) or simply $\delta$ ($\Delta$) the minimum (maximum) degree of a graph $G$. If $\delta = d_1 \leq \cdots \leq d_n = \Delta$ is a non-decreasing degree sequence of a graph $G$, then for convenience we will use the notation $(0^{x_0}, 1^{x_1}, \ldots, \Delta^{x_\Delta})$, where $x_k$ is the number of vertices of degree $k$ in the graph $G$. A cycle, a complete graph of order $n$, and a complete bipartite graph with cardinalities of partition sets $m$ and $n$ will be denoted by $C_n$, $K_n$, and $K_{m,n}$, respectively.

Let $G$ and $H$ be two disjoint graphs. We denote by $G + H$ the disjoint union of $G$ and $H$, which is a graph with vertex set $V(G) \cup V(H)$ and edge set $E(G) \cup E(H)$. Disjoint union of $k$ copies of graph $G$ is denoted by $kG$. We denote by $G \vee H$ the join of $G$ and $H$, which is a graph obtained from the disjoint union of $G$ and $H$ by adding edges joining every vertex of $G$ with every vertex of $H$. By $G + v$ we denote a graph $G$ with isolated vertex $v$, by $G + uv$ we denote a graph $G$ with the additional edge $uv$, connecting two nonadjacent vertices $u, v \in V(G)$, and by $G - uv$ we denote a graph $G$ without the edge $uv \in E(G)$. The number of components of $G$ is denoted by $c(G)$. We denote by $G[X]$, the subgraph $G$ induced by $X$. A subset $C$ of vertices of connected graph $G$ is called a vertex cutset of $G$ if $G - C := G[V \setminus C]$ is disconnected that is $c(G - C) \geq 2$.

Graph $G$ is called Hamiltonian if it has a simple cycle containing all vertices of the graph $G$. Graph $G$ is called pancyclic if it has simple cycles of all lengths from 3 to $n$. It is well-known, recognition problems, whether the graph is Hamiltonian or pancyclic, are NP-complete.

More than 40 years ago Chvátal introduced a new graph invariant, which he called graph toughness [12]. From then on a lot of research has been conducted, mainly related to the relationship between toughness conditions and the existence of cyclic structures [1–7, 10, 14–16, 18, 19].

Spectral conditions of pancyclicity for $t$-tough graphs

Historically, most of the research was based on a number of conjectures in [12]. The most challenging of these conjectures is still open: is there a finite constant $t_0$ such that every $t_0$-tough graph contains a Hamiltonian cycle? For a long time it was believed that this conjecture should hold for $t_0 = 2$. But in 2000, it was shown that the 2-tough conjecture is false [6]. On the other hand, we now know that the more general $t_0$-tough conjecture is true for a number of graph classes, including planar graphs, claw-free graphs, and chordal graphs. The early research in this area concentrated on sufficient degree conditions which, combined with a certain level of toughness, would yield the existence of long cycles.

Another stream involved finding toughness conditions for the existence of certain $k$-factors in graphs. For example, Enomoto, Jackson, and Sajto [14] have proved that each 2-tough graph of order at least 3 contains the 2-factor and for any $\varepsilon > 0$, there exists $(2 - \varepsilon)$-tough graph without the 2-factor and, hence, does not contain a Hamiltonian cycle.

The research on toughness has also focused on computational complexity issues. In particular, it has been shown that deciding whether $t(G) \geq t$ is coNP-complete for any fixed rational $t > 0$ [1].

## 2 Main Results

**Definition 2.1.** *For some real number $t$ graph $G$ is called $t$-tough if for any its vertex cutset $S$ the inequality $t \cdot c(G - S) \leq |S|$ is valid.*

**Definition 2.2.** *Toughness of graph $G$ is the greatest value $t$ for which $G$ is $t$-tough, and it is denoted by $\tau(G)$. If $G$ is the complete graph $K_n$ suppose $\tau(K_n) = +\infty$ for all $n \geq 1$.*

If $G$ is not complete graph, then toughness of the graph is equal to

$$\tau(G) = \min\left\{\frac{|S|}{c(G - S)} \ : \ S \subseteq V(G),\ c(G - S) \geq 2\right\},$$

where the minimum is taken over all vertex cutsets of $G$.

Obviously, $t$-tough graph is $s$-tough for any $s \leq t$ as well.

We need the following propositions:

**Proposition 2.1.** [13] *If $G \subset H$, where $G$ is a spanning subgraph, then $\tau(G) \leq \tau(H)$.*

**Proposition 2.2.** [13] *If $G$ is not complete graph, then $\tau(G) \leq \kappa(G)/2$, where $\kappa(G)$ is vertex connectivity of the graph.*

**Proposition 2.3.** [16] *If $G$ is not complete $t$-tough graph with a minimum degree $\delta(G)$, then $2t \leq \delta(G)$.*

**Proposition 2.4.** [17] *Let $G = (V(G), E(G))$ be a 2-tough graph of order $n$ and $x$ and $y$ be its two nonadjacent vertices, such that $d(x) + d(y) \geq n - 1$. Then $G' = G + xy$ is Hamiltonian if and only if $G$ is Hamiltonian.*

**Definition 2.3.** *For a nonnegative integer $k$ the $k$-closure of a graph $G$ of order $n$ is called the graph $\mathcal{C}_k(G)$ obtained from $G$ by recursively joining pairs of nonadjacent vertices the degree sum of which is at least $k$ as far as such a pair remains.*

Bondy and Chvatal showed that $\mathcal{C}_k(G)$ is uniquely defined by a given graph $G$ [9].

From Propositions 2.1 and 2.4 it follows that the Hamiltonian property is $(n-1)$-stable in the class of 2-tough graphs:

**Proposition 2.5.** *2-tough graph $G$ is Hamiltonian if and only if its $(n-1)$-closure $\mathcal{C}_{n-1}(G)$ is Hamiltonian.*

One more analogue of the Bondy-Chvatal closure theorem has been proved in [17] as well.

**Lemma 2.1.** ($t$-*closure lemma*) [17] *When $t \geq 2\left(\dfrac{3t-1}{2}\right)$-tough graph $G$ is Hamiltonian if and only if its $(n-t)$-closure $\mathcal{C}_{n-t}(G)$ is Hamiltonian.*

Let us denote the adjacency matrix of a graph $G$ by $A(G) = (a_{ij})$, where

$$a_{ij} = \begin{cases} 1, & \text{if } ij \in E(G), \\ 0, & \text{if } ij \notin E(G). \end{cases}$$

Let $D(G) = \text{diag}(d_G(v_1), \ldots, d_G(v_n))$ be the degree diagonal matrix of $G$. Then the matrix $Q(G) = D(G) + A(G)$ is called the signless Laplace matrix (Laplacian) of $G$. The set of all eigenvalues of matrix $A(G)$ is called the spectrum of $G$ and the greatest eigenvalue, denoted by $\mu(G)$, is called its spectral radius. The greatest eigenvalue of the signless Laplacian of $G$, denoted by $q(G)$, is called spectral radius of the signless Laplacian of $G$.

In recent years, the study of the Hamiltonian problem using the spectrum graph theory has received extensive attention, and some meaningful results have been obtained.

It is natural to assume that a $t$-tough graph is Hamiltonian or pancyclic if it possesses some additional properties. Recently, some sufficient conditions have been obtained for a $t$-tough graph to be Hamiltonian when $t \in \{1; 2; 3\}$.

**Proposition 2.6.** [11] *Let $G$ be a $t$-tough simple connected graph of order $n$ and size $m$, where $t \in \{1; 2; 3\}$. If $\mu(G) \geq \sqrt{n^2 - 4tn - 2n + 10t^2 + 2t - 1}$, then graph $G$ is Hamiltonian*

1. *when $t \in \{1; 2\}$ and $n \geq 8t$;*
2. *when $t = 3$ and $n > 9t$.*

In this work we establish sufficient conditions for a $t$-tough graph to be pancyclic in terms of the edge number, the spectral radius and the signless Laplacian spectral radius of the graph.

Let $G$ be a graph with degree sequence $d_1 \leq d_2 \leq \ldots \leq d_n$ and let $t$ be a positive integer. Let us formulate the following predicate $P(t)$ for degree sequence of $G$:

$$P(t): \ \forall i, t \leq i < \frac{n}{2}, \ d_i \leq i \Rightarrow d_{n-i+t} \geq n - i.$$

The following statements will be necessary for us as well.

**Proposition 2.7.** [16] *Let $t \in \{1; 2; 3\}$. If a $t$-tough graph satisfies $P(t)$, then it is Hamiltonian.*

Bondy observed that most conditions ensuring a graph $G$ is hamiltonian also ensures that $G$ is pancyclic or belongs to some easily described families of graphs (such as bipartite graphs) [19].

**Proposition 2.8.** [16] *Let $t$ be a positive integer. If a $t$-tough graph $G$ satisfies $P(t)$ and is Hamiltonian, then $G$ is pancyclic, or bipartite.*

**Proposition 2.9.** [16] *Let Hamiltonian graph $G$ of order $n$ contain more than $n/3$ vertices of degree greater than $n/2$. Then graph $G$ is pancyclic.*

**Proposition 2.10.** [16] *Let Hamiltonian graph $G$ of order $n$ contain at least $n/2$ vertices of degrees at least $n/2$. Then graph $G$ is pancyclic or bipartite or the graph $S_n$, where $S_n$ is the graph obtained from a clique $K_{n/2}$ and perfect matching $P$ on $n/2$ vertices by adding a perfect matching between $K_{n/2}$ and $P$ (see Fig. 1).*

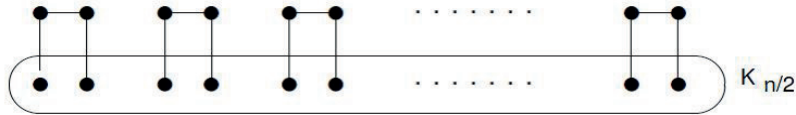

Figure 1: Graph $S_n$

The following result is of an independent interest.

**Theorem 2.1.** *Let $G$ be a simple connected $t$-tough graph of order $n$ and size $m$. If*

$$m \geq \binom{n - 2t}{2} + 3t^2,$$

*then graph $G$ is pancyclic or bipartite*

*1) when $t = 1$ and $n \geq 7$;*

*2) when $t = 2$ and $n \geq 16$;*

*3) when $t = 3$ and $n \geq 28$.*

*Proof.* Assume, to the contrary, that graph $G$ is neither pancyclic nor bipartite satisfying the conditions of the theorem. Then it does not satisfy the predicate $P(t)$, since otherwise, by Proposition 2.7, it would be Hamiltonian, and therefore, by Proposition 2.8, it would be pancyclic or bipartite. Thus, there is a positive integer $k$, such that $t \leq k < \frac{n}{2}$, $d_k \leq k$ and $d_{n-k+t} \leq n-k-1$. Then we have

$$2m = \sum_{i=1}^{k} d_i + \sum_{i=k+1}^{n-k+t} d_i + \sum_{i=n-k+t+1}^{n} d_i \leq$$
$$\leq k^2 + (n-2k+t)(n-k-1) + (k-t)(n-1) =$$
$$= n^2 - n + 3k^2 + k(1-2n-t) =$$
$$= 2\binom{n-2t}{2} + 6t^2 - (k-2t)(2n-3k-5t-1), \quad (1)$$

whence, by the theorem conditions, we obtain the inequalities:

$$\binom{n-2t}{2} + 3t^2 \leq m \leq \binom{n-2t}{2} + 3t^2 - \frac{(k-2t)(2n-3k-5t-1)}{2}, \quad (2)$$

thus $(k-2t)(2n-3k-5t-1) \leq 0$.

We consider the following three cases, according to the values of $t$.

**Case 1.** $t = 1$.

In this case $n \geq 7$ and $(k-2)(2n-3k-6) \leq 0$, then by Proposition 2.3 $\delta(G) \geq 2t = 2$, and therefore $2 \leq \delta(G) \leq d_k \leq k$. There are two subcases.

**Subcase 1.1.** If $(k-2)(2n-3k-6) = 0$, then $k = 2$, or $k \neq 2$ and $(2n-3k-6) = 0$.

**Subcase 1.1.1.** If $k = 2$, then $d_2 \leq 2$, $d_{n-1} \leq n-3$ and $d_n \leq n-1$. Hence, all inequalities in (2) which have the form $\binom{n-2}{2} + 3 \leq m \leq \binom{n-2}{2} + 3$, become equalities and graph $G$ has degree sequence: $(2^2, (n-3)^{n-3}, (n-1)^1)$. Then for graph $G$ there are two possibilities:

1) two vertices $w, z \in V(G)$ of degree 2 are nonadjacent. Then $G = K_1 \vee (K_{n-3} - uv + uw + vz)$ for some vertices $u, v \in V(K_{n-3})$, but it is easy to see that such a graph $G$ is Hamiltonian, and then by Proposition 2.9 $G$ is pancyclic graph, which is a contradiction.

2) two vertices $w, z \in V(G)$ of degree 2 are adjacent. Then $G = K_1 \vee (K_{n-3} + K_2)$, but it is easy to see that for such a graph $\tau(G) \leq \frac{1}{2}$, which is a contradiction.

**Subcase 1.1.2.** If $k \neq 2$ and $(2n-3k-6) = 0$, then from inequality $4k < 2n = 3k+6$ and the parity of the right side of this equality, we obtain $k = 4$, hence $n = 9$. Then $d_4 \leq 4, d_6 \leq 4, d_9 \leq 8$ and inequalities in (2) which have the form $48 \leq \sum_{i=1}^{9} d_i \leq 48$, become equalities, that is graph $G$ has degree sequence $(4^6, 8^3)$. This means that $G = K_3 \vee 3K_2$ which, it is easy to show, is Hamiltonian. Then by Proposition 2.9 $G$ is pancyclic graph, which is a contradiction.

**Subcase 1.2.** If $(k-2)(2n-3k-6) < 0$, then $k \geq 3$ and $(2n-3k-6) < 0$. Since $n \geq 2k+1$, then from inequalities $4k+2 \leq 2n < 3k+6$ and $k \geq 3$ it follows that $k = 3$. Therefore $n = 7$, and $m = 13$. Hence $d_3 \leq 3$, $d_5 \leq 3$, $d_7 \leq 6$. It is easy to verify that then there are the following graphical degree sequences:

$$(2^1, 3^4, 6^2), \; (3^5, 5^1, 6^1).$$

The first degree sequence corresponds to graph $G = K_2 \vee (2K_2 + K_1)$. But for such a graph $\tau(G) \leq \frac{2}{3}$, which is a contradiction. The second degree sequence corresponds to the graph $G$ represented in Fig. 2. It is easy to verify that such a graph is pancyclic, which is a contradiction.

**Case 2.** $t = 2$.

In this case $n \geq 16$ and $(k-4)(2n-3k-11) \leq 0$, then by Proposition 2.3 $\delta(G) \geq 2t = 4$, hence $4 \leq \delta(G) \leq d_k \leq k$. There are two subcases.

**Subcase 2.1.** If $(k-4)(2n-3k-11) = 0$, then $k = 4$, or $k \neq 4$ and $(2n-3k-11) = 0$.

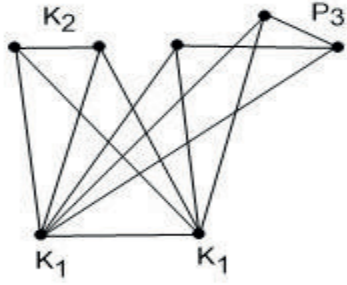

Figure 2: Graph with degree sequence $(3^5, 5^1, 6^1)$

**Subcase 2.1.1.** If $k = 4$, then $d_4 \leq 4$, $d_{n-2} \leq n - 5$ and $d_n \leq n - 1$. Hence, inequalities in (2) which have the form $\binom{n-4}{2} + 12 \leq m \leq \binom{n-4}{2} + 12$, are equalities and therefore, graph $G$ has degree sequence $(4^4, (n-5)^{n-6}, (n-1)^2)$.

Let us show that $(n-1)$-closure $\mathcal{C}_{n-1}(G)$ of graph $G$ is the complete graph $\mathcal{C}_{n-1}(G) = K_n$ when $n \geq 9$. Indeed, all $(n-6)$ vertices of degree $(n-5)$ must be pairwise adjacent in $\mathcal{C}_{n-1}(G)$. Therefore they have degrees at least $(n-5) + 2 = (n-3)$ in $\mathcal{C}_{n-1}(G)$. Then each of these vertices must be adjacent to each of four vertices of degree 4 in $\mathcal{C}_{n-1}(G)$. Hence, these four vertices have degrees at least $(n-6) + 2 = (n-4)$ in $\mathcal{C}_{n-1}(G)$. Therefore, each of these four vertices must be pairwise adjacent in $\mathcal{C}_{n-1}(G)$, that is $\mathcal{C}_{n-1}(G) = K_n$.

Then by Proposition 2.5 $G$ is Hamiltonian graph. Hence from $n > 10$ and by Proposition 2.9 graph $G$ is pancyclic, which is a contradiction.

**Subcase 2.1.2.** If $k \neq 4$ and $(2n - 3k - 11) = 0$ then it is easy to see that from inequality $4k < 2n = 3k + 11$ we obtain $k \leq 10$ and hence $16 \leq n \leq 20$. Then from the parity of the right part of the last equality it follows that either $k = 9$ and $n = 19$, or $k = 7$ and $n = 16$.

In case $k = 9$, $n = 19$ we have $d_9 \leq 9$, $d_{12} \leq 9$, $d_{19} \leq 18$. Then all inequalities in (2) which have the form $117 \leq m \leq 117$, are equalities and graph $G$ of order 19 has the degree sequence $(9^{12}, 18^7)$. By the similar arguments as above, it is easy to verify that 18-closure $\mathcal{C}_{18}(G)$ of the 2-tough graph $G$ is the complete graph $\mathcal{C}_{18}(G) = K_{19}$. Hence, by Proposition 2.5 $G$ is Hamiltonian graph. Then by Proposition 2.9 graph $G$ is pancyclic, which is a contradiction.

In case $k = 7$, $n = 16$ we have $d_7 \leq 7, d_{11} \leq 8, d_{16} \leq 15$. Hence, all inequalities in (2) which have the form $156 \leq m \leq 156$, are equalities and graph $G$ of order 16 has the degree sequence $(7^7, 8^4, 15^5)$.

Then it is easy to show that 15-closure $\mathcal{C}_{15}(G)$ of the 2-tough graph $G$ is the complete graph: $\mathcal{C}_{15}(G) = K_{16}$. Indeed, each of all four vertices of degree 8 must be pairwise adjacent and each of them must be adjacent to each of seven vertices of degree 7 in $\mathcal{C}_{15}(G)$. Therefore degrees of the last seven vertices must be at least 9. Then they must be pairwise adjacent in $\mathcal{C}_{15}(G)$. But it means that $\mathcal{C}_{15}(G) = K_{16}$.

Then by Proposition 2.5 $G$ is Hamiltonian graph. Since $G \neq S_{16}$, then by Proposition 2.10 graph $G$ is pancyclic or bipartite, which is a contradiction.

**Subcase 2.2.** If $(k-4)(2n - 3k - 11) < 0$, then $k \geq 5$ and $(2n - 3k - 11) < 0$. Then from inequalities $4k + 2 \leq 2n < 3k + 11$ it follows that $k \leq 8$. Therefore $16 \leq n \leq 17$, and $k = 8$, $n = 17$. Then $d_8 \leq 8, d_{11} \leq 8, d_{17} \leq 16$ and inequalities in (2), have the form $90 \leq m \leq 92$.

**Subcase 2.2.1.** In case $m = 92$ we have the unique degree sequence $(8^{11}, 16^6)$ to which graph $G = P \vee K_6$ corresponds, where $P$ is a 2-factor of order 11. By the similar arguments as above, it is easy to see that 16-closure $\mathcal{C}_{16}(G)$ of graph $G$ is the complete graph $\mathcal{C}_{16}(G) = K_{17}$. It means that, by Proposition 2.5 $G$ is Hamiltonian graph. Then by Proposition 2.9 graph $G$ is pancyclic, which is a contradiction.

**Subcase 2.2.2.** In case $m = 91$ it is easy to verify that there are the following graphical degree sequences:

$$(8^{11}, 15^2, 16^4), (7^2, 8^9, 16^6), (7^1, 8^{10}, 15^1, 16^5)$$

obtained from the graph $G = P \vee K_6$ removing one edge, and also degree sequences: $(6^1, 8^{10}, 16^6)$ corresponding to graph $G = (P + K_1) \vee K_6$, where $P$ is a 2-factor of order 10, and $(8^{11}, 14^1, 16^5)$, corresponding to graph $G$, shown in Figure 3.

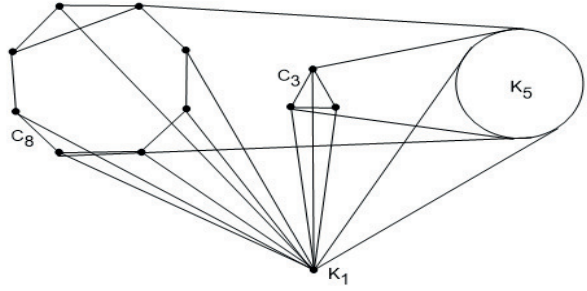

Figure 3: Graph with degree sequence $(8^{11}, 14^1, 16^5)$

Table 1: Graphical degree sequences.

| No | Degree sequence |
|---|---|
| 1 | $(5^1, 8^{10}, 15^1, 16^5)$ |
| 2 | $(6^1, 7^2, 8^8, 16^6)$ |
| 3 | $(6^1, 7^1, 8^9, 15^1, 16^5)$ |
| 4 | $(6^1, 8^{10}, 14^1, 16^5)$ |
| 5 | $(6^1, 8^{10}, 15^2, 16^4)$ |
| 6 | $(6^2, 8^9, 16^6)$ |
| 7 | $(7^4, 8^7, 16^6)$ |
| 8 | $(7^3, 8^8, 15^1, 16^5)$ |
| 9 | $(7^2, 8^9, 15^2, 16^4)$ |
| 10 | $(7^2, 8^9, 14^1, 16^5)$ |
| 11 | $(7^1, 8^{10}, 13^1, 16^5)$ |
| 12 | $(7^1, 8^{10}, 14^1, 15^1, 16^4)$ |
| 13 | $(7^1, 8^{10}, 15^3, 16^3)$ |
| 14 | $(8^{11}, 12^1, 16^5)$ |
| 15 | $(8^{11}, 14^2, 16^4)$ |

By the similar arguments as above, it is easy to show that 16-closures $\mathcal{C}_{16}(G)$ of all 2-tough graphs $G$ with the enumerated degree sequences coincide with the complete graph $\mathcal{C}_{16}(G) = K_{17}$. Hence, by Proposition 2.5 all these graphs $G$ are Hamiltonian. Then by Proposition 2.9 all these graphs $G$ are pancyclic, which is a contradiction.

**Subcase 2.2.3.** In case $m = 90$ it is easy to verify that there are the following degree sequences shown in Table 1, which are graphical.

By the similar arguments as above, it is easy to show that 16-closures $\mathcal{C}_{16}(G)$ of all 2-tough graphs, corresponding to these degree sequences, coincide with the complete graph $\mathcal{C}_{16}(G) = K_{17}$. Hence, by Proposition 2.5 all these graphs $G$ are Hamiltonian. Then by Proposition 2.9 all these graphs $G$ are pancyclic, which is a contradiction.

**Case 3.** $t = 3$.

In this case $(k-6)(2n-3k-16) \leq 0$. By Proposition 2.3 we have $6 \leq \delta(G) \leq d_k \leq k$. Besides, by the theorem condition the order of graph $G$ satisfies the inequality $n > 27$. There are two subcases.

**Subcase 3.1.** If $(k-6)(2n-3k-16) = 0$, then $k = 6$, or $k \neq 6$ and $2n - 3k - 16 = 0$.

**Subcase 3.1.1.** If $k = 6$, then $d_6 \leq 6$, $d_{n-3} \leq n-7$ and $d_n \leq n-1$. Hence, inequalities in (2) which have the form $\binom{n-6}{2} + 27 \leq m \leq \binom{n-6}{2} + 27$, are equalities and graph $G$ has degree sequence $(6^6, (n-7)^{n-9}, (n-1)^3)$.

Since graph $G$ is 3-tough then it is 2,5-tough. Hence, by Lemma 2.1 (when $t = 2$) it is Hamiltonian if and only if its $(n-2)$-closure $\mathcal{C}_{n-2}(G)$ is Hamiltonian. But it is easy to see that its $(n-2)$-closure is the complete graph: $\mathcal{C}_{n-2}(G) = K_n$. Thus, graph $G$ with degree sequence $(6^6, (n-7)^{n-9}, (n-1)^3)$ is Hamiltonian. Therefore, since $n > 14$, then by Proposition 2.9 graph $G$ is pancyclic, which is a contradiction.

**Subcase 3.1.2.** If $k \neq 6$ and $2n - 3k - 16 = 0$ then from the inequality $4k + 2 \leq 2n = 3k + 16$ it follows that $k \leq 14$, hence $28 \leq n \leq 29$. But, since $k$ and $n$ are integers, then from equality $2n - 3k - 16 = 0$ it follows that $n = 29$, $k = 14$. Then $d_{14} \leq 14, d_{18} \leq 14, d_{29} \leq 28$ and from inequalities (2) it follows that $m = 280$ and graph $G$ has degree sequence $(14^{18}, 28^{11})$. It is easy to show, that 27-closure of a 3-tough graph $G$ is the complete graph: $C_{27}(G) = K_{29}$. Hence, by Lemma 2.1 graph $G$ is Hamiltonian. Then by Proposition 2.9 graph $G$ is pancyclic, which is a contradiction.

**Subcase 3.2.** If $(k-6)(2n-3k-16) < 0$, then $k \geq 7$ and $2n - 3k - 16 < 0$. Then from inequalities $4k + 2 \leq 2n < 3k + 16$ it follows, that $k \leq 13$. Therefore $n \leq 27$, that contradicts the theorem conditions, and the proof is complete. □

**Lemma 2.2.** [14] *Let $G$ be a connected graph of order $n$ and size $m$. Then*

$$\rho(G) \leq \sqrt{2m - n + 1}.$$

*Equality occurs if and only if $G = K_n$ or $G = K_{1,n-1}$.*

**Lemma 2.3.** [13] *If $m \leq n$, then $\tau(K_{m,n}) = m/n$.*

Theorem 2.1 and Lemmas 2.2 and 2.3 imply the following proposition.

**Theorem 2.2.** *Let $G$ be a simple connected $t$-tough graph of order $n$ and size $m$, where $t \in \{1, 2, 3\}$. If the spectral radius of the graph satisfies the following inequality*

$$\rho(G) \geq \sqrt{n^2 - 4tn - 2n + 10t^2 + 2t + 1},$$

*then graph $G$ is pancyclic or bipartite,*

*1) when $t = 1$ and $n \geq 7$;*

*2) when $t = 2$ and $n \geq 16$;*

*3) when $t = 3$ and $n \geq 28$.*

*Proof.* By the theorem conditions and Lemma 2.3 the graph $G$ is different from graphs $K_n$ and $K_{1,n-1}$. Therefore we have inequalities

$$\sqrt{n^2 - 4tn - 2n + 10t^2 + 2t + 1} \leq \rho(G) < \sqrt{2m - n + 1},$$

whence we obtain

$$m > \binom{n - 2t}{2} + 3t^2.$$

Then by Theorem 2.1 the graph $G$ is pancyclic or bipartite. □

**Lemma 2.4.** [20] *Let $G$ be a connected graph of order $n$ and size $m$. Then the following inequality is valid:*

$$q(G) \leq \frac{2m}{n-1} + n - 2,$$

*and the equality holds if and only if either $G = K_n$ or $G = K_{1,n-1}$.*

**Theorem 2.3.** *Let $G$ be a simple $t$-tough connected graph of order $n$ and size $m$, where $t \in \{1, 2, 3\}$. If the spectral radius of the signless Laplacian of the graph satisfies the inequality*

$$q(G) \geq \frac{2n^2 + 10t^2 - 4tn + 2t - n}{n-1} + n - 2,$$

*then graph $G$ is pancyclic or bipartite,*

*1) when $t = 1$ and $n \geq 7$;*

*2) when $t = 2$ and $n \geq 16$;*

*3) when $t = 3$ and $n \geq 28$.*

*Proof.* By the theorem condition and Lemma 2.3 graph $G$ is different from graphs $K_n$ and $K_{1,n-1}$. Therefore we have inequalities

$$\frac{2n^2 + 10t^2 - 4tn + 2t - n}{n-1} + n - 2 \leq q(G) < \frac{2m}{n-1} + n - 2,$$

Spectral conditions of pancyclicity for $t$-tough graphs

whence we obtain
$$m > \binom{n-2t}{2} + 3t^2.$$

Then by Theorem 2.1 the graph $G$ is pancyclic or bipartite. □

## Acknowledgments

This work was supported by the Institute of Mathematics of the NAS of Belarus within the framework of the SPFR "Convergence – 2025".